\documentclass[10pt,a4paper,leqno]{amsart}

\date{}

\usepackage{latexsym}

\usepackage{amsfonts}

\usepackage{amsmath}

\usepackage{amssymb}

\usepackage{amscd}

\input xy
\xyoption{all}

\newcommand{\ma}[1]{\ensuremath{\mathbb{#1}}}
 
\font\bb=msbm7 at 10 pt

\def \C {\hbox{\bb C}}
\def \Z {\hbox{\bb Z}}
\def \Q {\hbox{\bb Q}}

\def \P {\mathcal{P}}
\def \F {\hbox{\bb F}}
\def \R {\hbox{\bb R}}
\def \A {\hbox{\bb A}}

\def \Tr {\mbox{\rm{Tr}}}

\def \T {\mathcal{T}}
\def \B {\mathcal{B}}

\def \L {\mathcal{L}}
\def \H {\mathcal{H}}
\def \U {\mathcal{U}}
\def \O {\mathcal{O}}
\def \K {\mathcal{K}}
\def \Q {\hbox{\bb Q}}
\def \I {\mbox{\bf I}}
\def \M {\mathcal{M}}
\def \II {\mathcal{I}}
\def \sgn {\mbox{\rm{sgn}}}

\newcommand{\Spec}{\ensuremath{\mbox{\rm{Spec }}}}

\newcommand{\Gal}{\ensuremath{\mbox{\rm{Gal }}}}

\newcommand{\Mat}{\ensuremath{\mbox{\rm{Mat}}}}

\begin{document}

\title{Newton stratification for polynomials: the open stratum.}



\author{R\'egis Blache}
\address{
\'Equipe ``G\'eom\'etrie Alg\'ebrique et Applications \`a la Th\'eorie de l'Information'', Universit\'e de Polyn\'esie Fran\c{c}aise, BP 6570, 98702 FAA'A, Tahiti, Polyn\'esie Fran\c{c}aise}
\email{blache@upf.pf}
\author{\'Eric F\'erard}
\address{
\'Equipe ``G\'eom\'etrie Alg\'ebrique et Applications \`a la Th\'eorie de l'Information'', Universit\'e de Polyn\'esie Fran\c{c}aise, BP 6570, 98702 FAA'A, Tahiti, Polyn\'esie Fran\c{c}aise}
\email{ferard@upf.pf}

\begin{abstract}

In this paper we consider the Newton polygons of $L$-functions 
coming from additive exponential sums associated to a polynomial over a finite field $\ma{F}_q$. These polygons define a stratification of the space of polynomials of fixed degree. We determine the open stratum: we give the generic Newton polygon for polynomials of degree $d\geq 2$ when the characteristic $p$ is greater than $3d$, and the Hasse polynomial, i.e. the equation defining the hypersurface complementary to the open stratum.
 
\end{abstract}

\maketitle

\section{Introduction}

Let $k:=\F_q$ be the finite field with $q:=p^m$ elements, and for any $r\geq 1$, let $k_r$ denote its extension of degree $r$. If  $\psi$ is a non trivial additive character on $\F_q$, then $\psi_r:=\psi\circ \Tr_{k_r/k}$ is a non trivial additive character of $k_r$, where $\Tr_{k_r/k}$ denotes the trace from $k_r$ to $k$. Let $f\in k[X]$ be a polynomial of degree $d\geq 2$ prime to $p$; then for any $r$ we form the additive exponential sum
$$S_r(f,\psi):=\sum_{x\in k_r}\psi_r(f(x)).$$
To this family of sums, one associates the $L$-function
$$L(f,T):=\exp\left(\sum_{r\geq 1} S_r(f,\psi)\frac{T^r}{r}\right).$$
It follows from the work of Weil on the Riemann hypothesis for function fields in characteristic $p$ that this $L$-function is actually a polynomial of degree $d-1$. Consequently we can write $$L(f,T)=(1-\theta_1T)\dots(1-\theta_{d-1}T).$$
Another consequence of the work of Weil is that the reciprocal roots $\theta_1,\dots,\theta_{d-1}$ are {\it $q$-Weil numbers of weight $1$}, i.e. algebraic integers all of whose conjugates have complex absolute $q^{\frac{1}{2}}$. Moreover, for any prime $\ell\neq p$, they are $\ell$-adic units, that is $|\theta_i|_\ell=1$. 

\medskip

A natural question is to determine their $q$-adic absolute value, or equivalently their $p$-adic valuation. In other words, one would like to determine the Newton polygon $NP_q(f)$ of $L(f,T)$ where $NP_q$ means the Newton polygon taken with respect to the valuation $v_q$ normalized by $v_q(q)=1$ ({\it cf.} \cite{ko}, Chapter IV for the link between the Newton polygon of a polynomial and the valuations of its roots). There is an elegant general answer to this problem when $p\equiv 1~[d]$, $p\geq 5$: then the Newton polygon $NP_q(f)$ has vertices ({\it cf.} \cite{ro}, Theorem 7.5) 
$$\left(n,\frac{n(n+1)}{2d}\right)_{1\leq n\leq d-1}.$$
This polygon is often called the {\it Hodge polygon} for polynomials of degree $d$, and denoted by $HP(d)$.

\medskip

Unfortunately, if we don't have $p\equiv 1~[d]$, there is no such general answer. We know that $NP_q(f)$ lies above $HP(d)$. This polygon can vary greatly depending on the coefficients of $f$, and it seems hopeless to give a general answer to the question above, as show the known examples ({\it cf.} \cite{sp} for degree $3$ polynomials, \cite{ho1} and \cite{ho2} for degree $4$ and degree $6$ polynomials respectively). On the other hand, we have asymptotic results ({\it cf.} \cite{zhu1}, \cite{zhu2}): in these papers, Zhu proves the one-dimensional case of Wan's conjecture ({\it cf.} \cite{wan} Conjecture 1.12), i.e. that there is a Zariski dense open subset $\U$ of the space of polynomials of degree $d$ over $\overline{\Q}$ such that when $p$ tends to infinity, for any $f\in \U$, the polygon $NP_q(f)$ obtained from the reduction of $f$ modulo a prime above $p$ in the field defined by the coefficients of $f$, we have $\lim_{p\rightarrow \infty} NP_q(f)=HP(d)$.

\medskip

A general result concerning Newton polygons is {\it Grothendieck's specialization theorem}. In order to quote it, let us recall some results about crystals. Let $\L_\psi$ denote the {\it Artin Schreier crystal}; this is an overconvergent $F$-isocrystal over $\A^1$ ({\it cf.} \cite{els} 6.5), and for any polynomial $f\in k[x]$ of degree $d$, we have an overconvergent $F$-isocrystal $f^*\L_\psi$ with ({\it cf.} \cite{bou})
$$L(f,T)=\det\left(1-T\phi_c|H^1_{\rm rig,c}(\A^1/K,f^*\L_\psi)\right).$$
Now if we parametrize the set of degree $d$ monic polynomials without constant coefficient by the affine space $\A^{d-1}$, associating the point $(a_1,\dots,a_{d-1})$ to the polynomial $f(X)=X^d+a_{d-1}X^{d-1}+\dots+a_1X$, we can consider the family of overconvergent $F$-isocrystals $f^*\L_\psi$. Now for a family of $F$-crystal $(\M,F)$ of rank $r$ over a $\F_p$-algebra $A$, we have Grothendieck's specialization theorem ({\it cf.} \cite{gr}, \cite{ka} Corollary 2.3.2)

\medskip

{\it Let $P$ be the graph of a continuous $\R$-valued function on $[0,r]$ which is linear between successive integers. The set of points in $\Spec(A)$ at which the Newton polygon of $(\M,F)$ lies above $P$ is Zariski closed, and is locally on $\Spec(A)$ the zero-set of a finitely generated ideal.}

\medskip

In other words, this theorem means that when $f$ runs over polynomials of degree $d$ over $\F_q$, then there is a Zariski dense open subset $U_{d,p}$ (the {\it open stratum}) of the (affine) space of these polynomials, and a {\it generic Newton polygon} $GNP(d,p)$ such that for any $f\in U_{d,p}$, $NP_q(f)=GNP(d,p)$, and $NP_q(f)\geq GNP(d,p)$ for any $f\in \F_q[X]$, $f$ monic of degree $d$ (where $NP\geq NP'$ means $NP$ lies above $NP'$). 

\medskip

The aim of this article is to determine explicitely both the generic polygon $GNP(d,p)$ and the associated {\it Hasse polynomial} $H_{d,p}$, i.e. the exact polynomial such that $U_{d,p}$ is the complementary of the hypersurface $H_{d,p}=0$. To be more precise, let $p\geq 3d$ be a prime; a {\it normalized} polynomial of degree $d$ over $\F_q$ is $f(x)=x^d+a_{d-2}x^{d-2}+\dots+a_1x\in \F_q[x]$; we identify the space of normalized polynomials with the affine space $\A^{d-2}(\F_q)$. Then the generic polygon $GNP(d,p)$ has vertices
$$\left(n,\frac{Y_n}{p-1}\right)_{1\leq n\leq d-1},~Y_n:=\min_{\sigma\in S_n} \sum_{k=1}^n \lceil\frac{pk-\sigma(k)}{d}\rceil,$$
and we have $NP_q(f)=GNP(d,p)$ exactly when $H_{d,p}(a_1,\dots,a_{d-2})\neq 0$, with $H_{d,p}$ the Hasse polynomial, that we determine explicitely. Note that both $GNP(d,p)$ and $H_{d,p}$ do not depend on $q$, but only on $p$. 

\medskip

The above results improve recent works of Scholten-Zhu ({\it cf.} \cite{sch}) and Zhu ({\it cf.} \cite{zhu1}, \cite{zhu2}). In \cite{sch}, Scholten and Zhu determine the first generic slope and the polynomials having this slope, and our work is a generalization of this result to the whole Newton polygon. In \cite{zhu1}, the generic Newton polygon is determined, but its $n$-th vertex depends on an intricated constant $\varepsilon_n$; moreover, Zhu doesn't need to give the exact equation defining $U_{d,p}$ since she just wants to prove its non emptyness. 

\medskip

We use $p$-adic cohomology, following the works of Dwork, Robba and others. To be more precise, we use Washnitzer-Monsky spaces of overconvergent series $\H^\dagger(A)$; one can define a linear operator $\beta$ on $\H^\dagger(A)$ and a differential operator $D$ with finite index on this space such that $\beta$ and $D$ commute up to a power of $p$. Then the linear map $\overline{\alpha}=\overline{\beta}^{\tau^{m-1}}\overline{\beta}^{\tau^{m-2}}\dots\overline{\beta}$ ($\tau$ being the Frobenius) on the quotient $\H^\dagger(A)/D\H^\dagger(A)$ has characteristic polynomial (almost) equal to $L(f,T)$. Using a monomial basis of $\H^\dagger(A)/D\H^\dagger(A)$, we are able to give congruences for the coefficients of the matrix $M:=\Mat_\B(\overline{\beta})$ in terms of the coefficients of a lift of $f$. We deduce congruences for the minors of $N:=\Mat_\B(\overline{\alpha})$, i.e. for the coefficients of the function $L(f,T)$.

\medskip

The paper is organized as follows: in section 1, we recall the results from $p$-adic cohomology we use, reducing the calculation of the $L$-function to the calculation of the matrix $N$. Section 2 is the technical heart of our work: we give congruences for the coefficients and the minors of $\Gamma$, a submatrix of $M$. Note that these results are sufficient to determine the generic Newton polygon in case $q=p$; moreover we deduce a congruence on exponential sums. In section $3$ we come to the general case: we give congruences for the minors of a submatrix $A$ of the matrix $N$, whose characteristic polynomial is $L(f,T)$. Finally we show the main results of the article in section 4, defining the generic Newton polygon for normalized polynomials of degree $d$ and the Hasse polynomial associated to this polygon ({\it cf.} Theorem 4.1).

\medskip 

\section{$p$-adic differential operators and exponential sums.}

In this section, we recall well known results about $p$-adic 
differential operators, and their application to the evaluation of the 
$L$-function of exponential sums. The reader interested in more details 
and the proofs should refer to \cite{ro}. 
 
\medskip

We denote by $\Q_p$ the field of $p$-adic numbers, and by $\K_m$ its 
(unique up to isomorphism) unramified extension of degree $m$. Let 
$\O_m$ be the valuation ring of $\K_m$; the elements of finite order in 
$\O_m^\times$ form a group $\T_m^\times$ of order $p^m-1$, and 
$\T_m:=\T_m^\times\cup\{0\}$ is the {\it Teichm\"uller} of $\K_m$. Note that 
it is the image of a section of reduction modulo $p$ from $\O_m$ to its 
residue field $\F_q$, called the {\it Teichm\"uller lift}. Let $\tau$ be the Frobenius; it is the 
generator of $\Gal(\K_m/\Q_p)$ which acts on $\T_m$ as the $p$th power 
map. Finally we denote by $\C_p$ a completion of a fixed algebraic 
closure $\overline{\Q}_p$ of $\Q_p$.

\medskip

Let $\pi \in \C_p$ be a root of the polynomial $X^{p-1}+p$. It is well known that $\Q_p(\pi)=\Q_p(\zeta_p)$ is a totally ramified extension of degree $p-1$ of $\Q_p$. We shall 
frequently use the valuation $v:=v_\pi$, normalized by $v_\pi(\pi)=1$, 
instead of the usual $p$-adic valuation $v_p$, or the $q$-adic valuation 
$v_q$.

\medskip

\subsection{Index of $p$-adic differential operators of order $1$.}

In this paragraph, we denote by $\Omega$ an algebraically closed field containing $\C_p$, complete under a valuation extending that of $\C_p$, and such that the residue class field of $\Omega$ is a transcendental extension of the residue class field of $\C_p$. For any $\omega \in \Omega$, $r\in \R$, we denote by $B(\omega,r^+)$ ({\it resp.} $B(\omega,r^-)$) the closed ({\it resp.} open) ball in $\Omega$ with center $\omega$ and radius $r$.

\medskip

Let $f(X):=\alpha_dX^d+\dots+\alpha_1X$, $\alpha_d\neq 0$ be a 
polynomial of degree $d$, prime to $p$, over the field $\F_q$, and let 
$g(x):=a_dX^d+\dots+a_1X \in \O_m[X]$ be the polynomial whose 
coefficients are the Teichm\"uller lifts of those of $f$.
Let $A:=B(0,1^+)\backslash B(0,1^-)$. We consider the space $\H^\dagger(A)$ of overconvergent 
analytic functions on $A$.

\medskip

Define the function $H:=\exp(\pi g(X))$; note that since $X\mapsto \exp(\pi X)$ has radius of convergence $1$, $H$ is not an element of $\H^\dagger(A)$. Now let $D$ be the differential 
operator (where a function acts on $\H^\dagger(A)$ by multiplication)
$$D:=X\frac{d}{dX}-\pi Xg'(X)~\left(=H^{-1}\circ 
X\frac{d}{dX}\circ H\right).$$ 
Since $H$ is not in $\H^\dagger(A)$, $D$ is injective in $\H^\dagger(A)$. Thus the index of $D$ in $\H^\dagger(A)$ is the dimension of its cokernel. By (\cite{ro} Proposition 
5.4.3 p226), this dimension is $d$.

\medskip

On the other hand, since $D$ can be seen as a differential operator acting 
on $\C_p[X,\frac{1}{X}]$, Theorem 5.6 of \cite{ro} ensures that a 
complementary subspace of $D\C_p[X,\frac{1}{X}]$ in 
$\C_p[X,\frac{1}{X}]$ is also a complementary subspace of 
$D\H^\dagger(A)$ in $\H^\dagger(A)$. Now an easy calculation gives, for 
any $n\in \Z$
$$DX^{n-d}=(n-d)X^{n-d}+\pi\sum_{i=1}^di\alpha_iX^{i+n-d},$$
and since this function is clearly in $D\H^\dagger(A)$, we get, for 
$n\geq d$
$$X^n\equiv -\frac{n-d}{\pi}X^{n-d}-\sum_{i=1}^{d-1} i\alpha_i 
X^{i+n-d}\quad [D\H^\dagger(A)],$$
and for $n<0$, $X^n\equiv -\frac{\pi}{n} \sum_{i=1}^di\alpha_iX^{i+n} 
~[D\H^\dagger(A)]$. Thus $\B:=\{1,\dots,X^{d-1}\}$ forms a basis of a 
complementary subspace of $D\H^\dagger(A)$ in $\H^\dagger(A)$, and for 
every $n\in \Z$, $X^n$ can be written uniquely as
$$X^n\equiv \sum_{i=0}^{d-1} a_{ni}X^i \quad [D\H^\dagger(A)],$$
for some $a_{ni}\in \K_m(\pi)$, $1\leq i\leq d-1$. We need more precise estimates for these coefficients and their 
$\pi$-adic valuations

\medskip

{\bf Lemma 1.1.} {\it We have the relations

i) $a_{ni}=\delta_{ni}$ if $0\leq n\leq d-1$,

ii) $v(a_{ni})\geq -\left[\frac{n-i}{d}\right]$ for $n\geq d$ and $i=1$

iii) $a_{n0}=0$ for any $n>0$.}

\medskip

{\it Proof.} Part {\it i)} is trivial, and part {\it ii)} is just Lemma 
7.7 in \cite{ro}. It remains to show part {\it iii)}; from the 
discussion above the lemma and the definition of the $a_{ni}$, we get 
for any $n\geq d$
$$a_{n0}= -\frac{n-d}{\pi}a_{n-d,0}-\sum_{i=1}^{d-1} i\alpha_i 
a_{i+n-d,0}.$$
Thus $a_{d0}=0$ from part {\it i)}, and the result follows recursively.

\subsection{L-functions of exponential sums as characteristic polynomials.}

We define the power series $\theta(X):=\exp(\pi X-\pi X^p)$; this is a {\it splitting function} in Dwork's terminology ({\it cf.} \cite{dw} p55). Its values at the points of $\T_1$ are $p$-th roots of unity; in other words this function represents an additive character of order $p$. It is well known that 
$\theta$ converges for any $x$ in $\C_p$ such that 
$v_p(x)>-\frac{p-1}{p^2}$, and in particular $\theta \in \H^\dagger(A)$. 
We will need the following informations on the coefficients of 
the power series $\theta$

\medskip

{\bf Lemma 1.2.} {\it Set $\theta(X):= \sum_{i\geq 0} b_iX^i$; then we have

i) $b_i=\frac{\pi^i}{i!}$ if $0\leq i\leq p-1$;

ii) $v(b_i)\geq i$ for $0\leq i\leq p^2-1$;

iii) $v(b_i)\geq \left(\frac{p-1}{p}\right)^2i$ for $i\geq p^2$.}

\medskip

We define the functions $F(X):=\prod_{i=1}^d 
\theta(a_iX^i):=\sum_{n\geq0} h_nX^n$, and $G(X):=\prod_{i=0}^{m-1} 
F^{\tau^i}(X^{p^i})$; since $\theta$ is overconvergent, $F$ and $G$ 
also, and we get $G\in \H^\dagger(A)$. 

\medskip

Consider the mapping $\psi_q$ defined on $\H^\dagger(A)$ by 
$\psi_qf(x):=\frac{1}{q}\sum_{z^q=x}f(z)$; if $f(X)=\sum b_nX^n$, then 
$\psi_q f(X)=\sum b_{qn}X^n$. Let $\alpha:=\psi_q \circ G$; as operators 
on $\H^\dagger(A)$, $D$ and $\alpha$ commute up to a factor $q$, and we 
get a commutative diagram with exact rows 
$$\xymatrix{
 0 \ar[r]& \ar[d]_{q\alpha} \H^\dagger(A) \ar[r]^{D} & \ar[d]_{\alpha} 
\H^\dagger(A) \ar[r] & \ar[d]_{\overline{\alpha}} \H^\dagger(A)/D\H^\dagger(A) \ar[r] & 0\\
0 \ar[r]& \H^\dagger(A) \ar[r]^{D} & \H^\dagger(A) \ar[r] & \H^\dagger(A)/D\H^\dagger(A) \ar[r] & 0\\
}$$
Let $L^*(f,T)$ be the $L$-function associated to the sums 
$S_r^*(f):=\sum_{x\in k_r^\times} \psi_r(f(x))$; Dwork's trace 
formula (cf \cite{ro}) gives the following
$$L^*(f,T)=\frac{\det(1-T\alpha)}{\det(1-qT\alpha)}=\det(1-T\overline{\alpha}).$$
We have thus rewritten the $L$-function associated to the family of 
exponential sums as the characteristic polynomial of an endomorphism in 
a $p$-adic vector space.

\medskip

Let $\beta$ be the endomorphism of $\H^\dagger(A)$ defined by 
$\beta=\psi_p\circ F$; then $\tau^{-1}\circ\beta$ commutes with $D$ up to a factor 
$p$, and passes to the quotient, giving an endomorphism 
$\overline{\tau^{-1}\circ\beta}$ of $W$, the $\K_m(\zeta_p)$-vector space with basis $\B$. Thus $\beta$ induces $\overline{\beta}$ from $W$ to $W^\tau$, the $\K_m(\zeta_p)$-vector space $W$ with scalar multiplication given by $\lambda\cdot w=\lambda^\tau w$. On the other hand we have $\alpha=\beta^{\tau^{m-1}}\dots \beta^{\tau}\beta$. This gives the following relation between the endomorphism  $\overline{\alpha}$ of $W$ and the semilinear morphism $\overline{\beta}$ (note that $W^{\tau^m}=W$)
$$\overline{\alpha}=\overline{\beta}^{\tau^{m-1}}\dots 
\overline{\beta}^{\tau}\overline{\beta}.$$

\medskip

Let $M:=Mat_\B(\overline{\beta})$ ({\it resp.} $N$) be the matrix of 
$\overline{\beta}$ ({\it resp.} $\overline{\alpha}$) in the basis $\B$, 
and $m_{ij}$ ({\it resp.} $n_{ij}$), $0\leq i,j\leq d-1$ be the 
coefficients of $M$ ({\it resp.} $N$). From the description of $F$, we 
can write $m_{ij}=h_{pi-j}+\sum_{n\geq d} h_{np-j}a_{ni}$ (cf \cite{ro} 
7.10). Since we have $h_0=1$, and $h_n=0$ for negative $n$, we see from 
Lemma 1.2 {\it iii)} that $m_{00}=1$, and $m_{0j}=0$ for $1\leq j\leq 
d-1$. Since $N=M^{\tau^{m-1}}\dots M$, the same is true for the $n_{0i}$; thus 
the space $W'=Vect(X,\dots,X^{d-1})$ is stable under the action of $\overline{\alpha}$,
({\it resp.} $\overline{\beta}$ induces a morphism from $W'$ to $W'^\tau$) and the matrix $\Gamma$ ({\it resp.} $A$) defined by $\Gamma:=\left(m_{ij}\right)_{1\leq i,j\leq d-1}$, 
({\it resp.} $A:=\left(n_{ij}\right)_{1\leq i,j\leq d-1}$) is the matrix of the restriction of $\overline{\beta}$ ({\it resp.} $\overline{\alpha}$) with respect to the basis $\{X,\dots,X^{d-1}\}$. These matrices satisfy
$A=\Gamma^{\tau^{m-1}}\dots \Gamma$, and 
$\det(1-T\overline{\alpha})=(1-T)\det(\I_{d-1}-TA)=(1-T)\det(\I_{d-1}-T\Gamma^{\tau^{m-1}}\dots 
\Gamma)$. Finally, since we assumed $f(0)=0$, we have 
$S_r^*(f)=S_r(f)-1$ for any $r\geq 1$, and $L^*(f,T)=(1-T)L(f,T)$. From this we deduce the 
following result, which we will use to evaluate the valuations of the 
coefficients of the $L$-function associated to $f$

\medskip

{\bf Proposition 1.1.} {\it Let $\Gamma$ be as above; then we have
$$L(f,T)=\det(\I_{d-1}-T\Gamma^{\tau^{m-1}}\dots \Gamma).$$}

\medskip

{\bf Remark 1.1.} We have chosen to work over a ring of 
overconvergent series, the Washnitzer-Monsky dagger space; one can check 
that if $K:=\K_m(\gamma)$ is the totally ramified extension of $\K_m$ containing a fixed root of $X^d-\pi$, then the space $W'\otimes K$ with $W'$ as above is isomorphic 
to the space $H_0(SK_{\bullet}(B,D))$ constructed in \cite{as}, and under this isomorphism
the operator $\overline{\alpha}$ corresponds to $H_0(\alpha)$ there. Moreover,
these spaces are isomorphic to the first rigid cohomology group 
$H^1_{\rm rig,c}(\A^1/K,f^*\L_\psi)$ ({\it cf.} \cite{bou}).
\bigskip

\section{Congruences for the coefficients and the minors of the matrix 
$\Gamma$.}

In this section, we express the ``principal parts" of the coefficients 
$m_{ij}$ in terms of certain coefficients of the powers of the lifting 
$g$ of the polynomial $f$. Then we use these results to give the 
principal parts of the coefficients of the $L$-function.

\subsection{The coefficients.} Recall that we can express the 
coefficients $m_{ij}$ from the coefficients $h_n$ of the power series 
$F$ and the $a_{ni}$ in the following way
$$m_{ij}=h_{pi-j}+\sum_{n\geq d} h_{np-j}a_{ni}.$$
We begin by a congruence on the coefficients of $F$.

\medskip

{\bf Notation.} Let $P$ be a polynomial; we denote by 
$\left\{P\right\}_n$ its coefficient of degree $n$.

\medskip

{\bf Lemma 2.1} {\it Assume $p\geq d$, and let $0\leq n\leq (p-1)d$; 
then we have the following congruence for the coefficients of the power 
series $F$
$$h_n\equiv  
\sum_{k=\lceil\frac{n}{d}\rceil}^{p-1}\left\{g^{k}\right\}_n\frac{\pi^{k}}{k!}\quad 
[p\pi],$$
where $\lceil r\rceil$ is the least integer greater or equal than $r$.}

\medskip

{\it Proof.} From the definition of $F$, we get
$$h_n=\sum_{m_1+\dots+dm_d=n}a_1^{m_1}\dots a_d^{m_d}b_{m_1}\dots b_{m_d}.$$
Since $m_1+\dots+dm_d=n$, we get $d(m_1+\dots+m_d)\geq n$, and 
$m_1+\dots+m_d\geq \lceil\frac{n}{d}\rceil$; on the other hand we 
clearly have $m_1+\dots+m_d\leq n$, and we write
$$h_n=\sum_{k=\lceil\frac{n}{d}\rceil}^n h_{n,k},\qquad 
h_{n,k}=\sum_{m_1+\dots+dm_d=n\atop{m_1+\dots+m_d=k}}a_1^{m_1}\dots 
a_d^{m_d}b_{m_1}\dots b_{m_d}.$$
From Lemma 1.2 {\it ii)}, since $n<pd\leq p^2$, we have $m_i<p^2$, and $v(b_{m_i})\geq m_i$; thus $v(h_{n,k})\geq 
k$, and $h_n\equiv \sum_{k=\lceil\frac{n}{d}\rceil}^{p-1} 
h_{n,k}~\left[p\pi\right]$. Since $k\leq p-1$, the same is true for the 
$m_i$ appearing in the expression of $h_{n,k}$: from Lemma 1.2 {\it i)}, 
we know the $b_{m_i}$ explicitely, and we get
$$h_{n,k}  =  
\sum_{m_1+\dots+dm_d=n\atop{m_1+\dots+m_d=k}}\frac{a_1^{m_1}\dots 
a_d^{m_d}\pi^{k}}{m_1!\dots m_d!} =  
\frac{\pi^{k}}{k!}\sum_{m_1+\dots+dm_d=n\atop{m_1+\dots+m_d=k}}\binom{k}{m_1,\dots,m_d}a_1^{m_1}\dots
a_d^{m_d}$$
where $\binom{k}{m_1,\dots,m_d}:=\frac{k!}{m_1!\dots m_d!}$ denotes a 
multinomial coefficient. On the other hand, developing the polynomial 
$g^{k}$ yields
$$g^{k}(X)=\left(\sum_{i=1}^d a_iX^i\right)^{k}=\sum_{m_1+\dots+m_d=k} 
\binom{k}{m_1,\dots,m_d}a_1^{m_1}\dots a_d^{m_d}X^{\sum im_i},$$
and we get the result.

\medskip

We now give a congruence on the coefficients $m_{ij}$ of $\Gamma$.

\medskip

{\bf Proposition 2.1} Assume that $p\geq d+3$. Let $1\leq i,j\leq d-1$; 
we have
$$m_{ij}\equiv h_{pi-j}~[p\pi].$$

\medskip

{\it Proof.} From the expression of $m_{ij}$, we are reduced to show 
that for any $n\geq d$, we have
$v(h_{np-j}a_{ni})\geq p$.

\medskip

Assume first that $n\leq p$; from the expression of $h_n$, we see that 
the $m_i$ appearing in $h_{np-j}$ are all less than $p^2-1$, and we have 
$v(h_{np-j})\geq \frac{np-j}{d}$. Let $d\leq n< d+i$; from Lemma 1.1, we 
have $v(a_{ni})\geq -\left[\frac{n-i}{d}\right]\geq 0$, and 
$v(h_{np-j}a_{ni})\geq\frac{np-j}{d}\geq \frac{dp-j}{d}>p-1$. On the 
other hand, if $n\geq d+i$, $v(a_{ni})\geq -\left[\frac{n-i}{d}\right]\geq 
\frac{i-n}{d}$, and 
$v(h_{np-j}a_{ni})\geq\frac{np-j}{d}+\frac{i-n}{d}=\frac{n(p-1)+i-j}{d}\geq 
p-1 +\frac{i(p-1)+i-j}{d}>p-1$ since $p\geq d$.

\medskip

Suppose now that $n>p$; in this case we have $v(h_{np-j})\geq 
\frac{np-j}{d}\left(\frac{p-1}{p}\right)^2$ (cf \cite{ro} Lemma on 
p242). Thus
$$v(h_{np-j}a_{ni})\geq 
\frac{np-j}{d}\left(\frac{p-1}{p}\right)^2-\frac{n-i}{d}=\frac{n}{d}\left(\frac{(p-1)^2}{p}-1\right)-\frac{1}{d}\left(\left(\frac{p-1}{p}\right)^2j-i\right).$$
We have $\left(\frac{p-1}{p}\right)^2j-i\leq d$, thus 
$v(h_{np-j}a_{ni})\geq \frac{n}{d}\left(\frac{(p-1)^2}{p}-1\right)-1$. 
Since $n>p$, we get $\frac{n}{p}>1$ and 
$v(h_{np-j}a_{ni})>\frac{p^2-3p+1}{d}-1>p-1$ for $p\geq d+3$.

\medskip

{\bf Corollary 2.1} Assume that $p\geq d+3$. Let $1\leq i,j\leq d-1$; we 
have
$$m_{ij}\equiv\left\{g^{\lceil\frac{pi-j}{d}\rceil}\right\}_{pi-j}\frac{\pi^{\lceil\frac{pi-j}{d}\rceil}}{\lceil\frac{pi-j}{d}\rceil!}\quad
\left[\pi^{\lceil\frac{pi-j}{d}\rceil+1}\right].$$

\medskip

Another consequence of the above evaluations is a congruence on 
exponential sums associated to polynomials over the prime field: since $S_1(f)$ is the trace of the matrix $\Gamma$, we 
deduce from proposition 2.1

\medskip

{\bf Corollary 2.2} {\it Assume $p\geq d+3$, and let $f\in \F_p[X]$ be a polynomial of degree $d$; then we have the following 
congruence on the exponential sum $S_1(f)$

$$S_1(f)\equiv 
\sum_{k=\lceil\frac{p-1}{d}\rceil}^{p-1}\sum_{i=1}^{d-1}\left\{g^k\right\}_{(p-1)i}\frac{\pi^k}{k!}~[p\pi].$$}

\medskip

\subsection{The minors.} Our aim here is to give estimates for the 
principal parts of certain minors of the matrix $\Gamma$. Recall 
the following expression of a characteristic polynomial
$$\det(\I_{d-1}-T\Gamma )=1+\sum_{n=1}^{d-1} M_nT^n,$$
where $M_n=\sum_{1\leq u_1<\dots<u_n\leq d-1} \sum_{\sigma\in S_n} 
\sgn(\sigma) \prod_{i=1}^n m_{u_iu_{\sigma(i)}}$ is the sum of the 
$n\times n$ minors centered on the diagonal of $\Gamma$. We use the 
results of paragraph 2.1 to give a congruence for the coefficients $M_n$.

\medskip

{\bf Definition 2.1} {\it {\it i)} Set $Y_n:= \min_{\sigma\in S_n} 
\sum_{k=1}^n \lceil\frac{pk-\sigma(k)}{d}\rceil$, and
$$\Sigma_n:=\{\sigma\in S_n,~\sum_{k=1}^n 
\lceil\frac{pk-\sigma(k)}{d}\rceil=Y_n\}.$$

{\it ii)} For every $1\leq i \leq d-1$, set $j_i$ be the least positive 
integer congruent to $pi$ modulo $d$, and for every $1\leq n \leq d-1$, 
let $B_n:=\{1\leq i\leq n,~j_i\leq n\}$.}

\medskip

Note that since $p$ is coprime to $d$, the map $i\mapsto j_i$ is an 
element of $S_{d-1}$, the $d-1$-th symetric group.  We can use the set 
$B_n$ to describe $\Sigma_n$ precisely

\medskip

{\bf Lemma 2.2.} {\it Let $1\leq n\leq d-1$; we have $\Sigma_n=\{ 
\sigma\in S_n,~\sigma(i)\geq j_i~\forall i\in B_n\}$, and 
$Y_n=\sum_{k=1}^n \lceil\frac{pk}{d}\rceil-\#B_n$.}

\medskip

{\it Proof.} It is easily seen that for any $1\leq j\leq j_i-1$, we have 
$\lceil\frac{pi-j}{d}\rceil=\lceil\frac{pi}{d}\rceil$, and for $j_i\leq 
j\leq n$, $\lceil\frac{pi-j}{d}\rceil=\lceil\frac{pi}{d}\rceil-1$. From 
this we deduce
$$\sum_{k=1}^n \lceil\frac{pk-\sigma(k)}{d}\rceil=\sum_{k=1}^n 
\lceil\frac{pk}{d}\rceil-\#\{1\leq k\leq n,~\sigma(k)\geq j_k\}.$$
Now we have the inclusion $\{1\leq k\leq n,~\sigma(k)\geq j_k\}\subset 
B_n$. Finally the set $\{ \sigma\in S_n,~\sigma(i)\geq j_i~\forall i\in 
B_n\}$ is not empty, since $i\mapsto j_i$ is an injection from $B_n$ 
into $\{1,\dots,n\}$; we get $Y_n=\sum_{k=1}^n 
\lceil\frac{pk}{d}\rceil-\#B_n$, and that the permutations reaching this 
minimum  are exactly the ones with $\sigma(i)\geq j_i$ for all $i\in 
B_n$. This is the desired result.

\medskip

We are now ready to give a congruence for the coefficients $M_n$ of the 
polynomial $\det(\I_{d-1}-T\Gamma)$.

\medskip

{\bf Definition 2.2.} {\it Recall that we have set $g(X)=\sum_{i=1}^d 
a_iX^i$. For any $1\leq n\leq d-1$ let $\P_n$ be the polynomial in 
$\Z[X_1,\dots,X_d]$ defined by
$$\P_n(a_1,\dots,a_d):=\sum_{\sigma\in \Sigma_n} \sgn(\sigma)
\prod_{i=1}^n\left\{g^{\lceil\frac{pi-\sigma(i)}{d}\rceil}\right\}_{pi-\sigma(i)}.$$}

\medskip

{\bf Lemma 2.3} {\it Let $1\leq u_1<\dots<u_n= n+s$ and $1\leq v_1<\dots<v_n= n+t$ be integers; then we have the following inequality
$$\sum_{k=1}^n \lceil\frac{pu_k-v_k}{d}\rceil\leq Y_n+\left(\left[\frac{p}{d}\right]-1\right)t-s.$$}

\medskip

{\it Proof.} We first rewrite the sum as in the proof of lemma 2.2
$$\sum_{k=1}^n \lceil\frac{pu_k-v_k}{d}\rceil=\sum_{k=1}^n \lceil\frac{pu_k}{d}\rceil-\#\{v_i,~v_i\geq j_{u_i}\}.$$
We know that there are $\#B_n$ integers in $\{1,\dots,n\}$ such that $j_i\leq n$ ; thus there are at most $\#B_n+s$ integers in $\{1,\dots,n\}$ such that $j_i\leq n+s$ since $i\mapsto j_i$ is a bijection. On the other hand, there are at most $t$ elements in $\{n+1,\dots,n+t\}$ such that $j_i\leq n+s$; thus the set $\#\{v_i,~v_i\geq j_{u_i}\}$ contains at most $\#B_n+s+t$ elements, and we get
$$\begin{array}{rcl}
\sum_{k=1}^n \lceil\frac{pu_k-v_k}{d}\rceil & \geq & 
\sum_{k=1}^n \lceil\frac{pu_k}{d}\rceil-\#B_n-s-t\\
& \geq & \sum_{k=1}^{n} 
\lceil\frac{pk}{d}\rceil+\lceil\frac{p(n+t)}{d}\rceil-\lceil\frac{pn}{d}\rceil-\# B_n -s-t\\
& \geq & Y_n+\lceil\frac{p(n+t)}{d}\rceil-\lceil\frac{pn}{d}\rceil-s-t.\\
\end{array}$$
Now for any $a,b\geq 0$ we have $\lceil a+b\rceil\geq \lceil 
a\rceil+[b]$, and the sum above is greater than 
$Y_n+\left[\frac{pt}{d}\right]-t-s$. Moreover, $[ab]\geq [a][b]$, and the 
sum is greater than $Y_n+\left(\left[\frac{p}{d}\right]-1\right)t-s$. This proves the lemma.

\medskip

{\bf Proposition 2.2} {\it Assume $p\geq 3d$; then for any $1\leq n\leq 
d-1$, we have
$$M_n\equiv \frac{\P_n(a_1,\dots,a_d)}{\prod_{i\notin 
B_n}\lceil\frac{pi}{d}\rceil!\prod_{i\in 
B_n}\left(\lceil\frac{pi}{d}\rceil-1\right)!}\pi^{Y_n}\quad[\pi^{Y_n+1}].$$}

\medskip

{\it Proof.} We first choose a term in the development of $M_n$ with 
$\{u_1,\dots,u_n\}\neq\{1,\dots,n\}$; let $u_n=n+t$, $t\geq 1$. From 
Corollary 2.1, we have
$$v(\prod_{k=1}^n m_{u_ku_{\sigma(k)}})\geq \sum_{k=1}^n 
\lceil\frac{pu_k-u_{\sigma(k)}}{d}\rceil.$$
Applying Lemma 2.3 to the $u_i$ and $v_i:=u_{\sigma(i)}$, we get that the 
valuation is greater than $Y_n+\left(\left[\frac{p}{d}\right]-2\right)t$. 
Finally since $p\geq 3d$ and $t\geq 1$, the valuation of the term above 
is greater than $Y_n+1$ and we need only consider the terms in the 
development of $M_n$ with $u_1=1,\dots,u_n=n$ to get the result.

\medskip

 From Corollary 2.1 and the description of $M_n$, we get
$$M_n\equiv \sum_{\sigma\in S_n} \sgn(\sigma)
\prod_{i=1}^n\left\{\frac{g^{\lceil\frac{pi-\sigma(i)}{d}\rceil}}{\lceil\frac{pi-\sigma(i)}{d}\rceil!}\right\}_{pi-\sigma(i)}\pi^{\sum_{i=1}^n\lceil\frac{pi-\sigma(i)}{d}\rceil}\quad[\pi^{Y_n+1}],$$
and we can restrict the sum to $\Sigma_n$ from the definition of $Y_n$. 
Finally for any $\sigma\in \Sigma_n$, we have 
$\lceil\frac{pi-\sigma(i)}{d}\rceil=\lceil\frac{pi}{d}\rceil$ if 
$i\notin B_n$, and 
$\lceil\frac{pi-\sigma(i)}{d}\rceil=\lceil\frac{pi}{d}\rceil-1$ else; 
thus the product $\prod_{i=1}^n \lceil\frac{pi-\sigma(i)}{d}\rceil!$ is 
independent of the choice of $\sigma$ in $\Sigma_n$. This ends the proof of 
Proposition 2.2.

\medskip

\section{Congruences for the minors of $A$.}

In this section, we give congruences for the coefficients of the 
characteristic polynomial of:
$$A=\Gamma^{\tau^{m-1}}\Gamma^{\tau^{m-2}}\dots\Gamma.$$
Recall $A:=(n_{ij})_{1\leq i,j\leq d-1}$, and set 
$\det(\I_{d-1}-TA):=\sum_{n=0}^{d-1} \M_nT^n$, with:
$$\M_n=\sum_{1\leq u_1<\dots<u_n\leq d-1} \sum_{\sigma\in 
S_n}\sgn(\sigma)\prod_{i=1}^n n_{u_i,u_{\sigma(i)}}.$$
Let us give an expression for $n_{ij}$:
$$n_{ij}=\sum_{1\leq k_1,\dots,k_{m-1}\leq 
d-1}m_{ik_1}^{\tau^{m-1}}m_{k_1k_2}^{\tau^{m-2}}\dots m_{k_{m-1}j}.$$
Fix $U=\{u_1,\dots,u_n\}$; replacing the above in $S_U:=\sum_{\sigma\in 
S_n}\sgn(\sigma)\prod_{i=1}^n n_{u_i,u_{\sigma(i)}}$, we get (where the inner
sum in the first line taken over $1\leq j\leq n$, and the other ones over $1\leq i\leq m-1$, $1\leq j\leq n$):
$$ S_U  =  \sum_{\sigma\in S_n}\sgn(\sigma)\prod_{i=1}^n \sum_{1\leq 
k_{ij}\leq 
d-1}m_{u_ik_{1i}}^{\tau^{m-1}}m_{k_{1i}k_{2i}}^{\tau^{m-2}}\dots 
m_{k_{m-1i}u_{\sigma(i)}} $$
$$\qquad =  \sum_{1\leq k_{ij}\leq d-1}\sum_{\sigma\in 
S_n}\sgn(\sigma)\prod_{i=1}^n 
m_{u_ik_{1i}}^{\tau^{m-1}}m_{k_{1i}k_{2i}}^{\tau^{m-2}}\dots 
m_{k_{m-1i}u_{\sigma(i)}}$$
$$\qquad =  \sum_{1\leq k_{ij}\leq d-1}\prod_{i=1}^n 
m_{u_ik_{1i}}^{\tau^{m-1}}\dots 
m_{k_{m-2i}k_{m-1i}}^{\tau}\sum_{\sigma\in S_n}\sgn(\sigma)\prod_{i=1}^n 
m_{k_{m-1i}u_{\sigma(i)}}$$

\medskip

{\bf Lemma 3.1} {\it If the map $i\mapsto k_{m-1i}$ is not injective, we 
have:
$$S':=\sum_{\sigma\in S_n}\sgn(\sigma)\prod_{i=1}^n 
m_{k_{m-1i}u_{\sigma(i)}}=0.$$}

\medskip

{\it Proof.} Assume that $k_{m-1i}=k_{m-1j}$ for some $i\neq j$. Then 
$\sigma\mapsto \sigma'=\sigma\circ(i,j)$ is a bijection from $A_n$ to 
$S_n\backslash A_n$, and we write
$$S'=\sum_{\sigma\in A_n} \left(\sgn(\sigma)\prod_{l=1}^n 
m_{k_{m-1l}u_{\sigma(l)}}+\sgn(\sigma')\prod_{l=1}^n 
m_{k_{m-1l}u_{\sigma'(l)}}\right)~;$$
Since $\sgn(\sigma')=-\sgn(\sigma)$, the sum above is zero for any $\sigma$.

\medskip

Thus we can write $k_{m-1i}=\theta_{m-1}(i)$ for some injective map 
$\theta_{m-1}:\{1,\dots,n\}\rightarrow\{1,\dots,d-1\}$. Let $\II_n$ be 
the set of such maps. We get a new expression for $S$ (where the first 
sum is taken over $1\leq i\leq m-2$, $1\leq j\leq n$)
$$S_U=\sum_{1\leq k_{ij}\leq d-1}\sum_{\theta_{m-1}\in 
\II_n}\sum_{\sigma\in S_n}\sgn(\sigma)\prod_{i=1}^n 
m_{u_ik_{1i}}^{\tau^{m-1}}m_{k_{1i}k_{2i}}^{\tau^{m-2}}\dots 
m_{\theta_{m-1}(i)u_{\sigma(i)}},$$

\medskip

Now we show that each of the maps $\theta_j:i\mapsto k_{ji}$ must be in 
$\II_n$:

\medskip

{\bf Lemma 3.2} {\it Assume that the maps $\theta_l:i\mapsto k_{li}$ are 
in $\II_n$ for any $1\leq t<l\leq m-1$, but that the map $i\mapsto k_{ti}$ is 
not injective; then we have the equality:
$$S'':=\sum_{(\theta_{t+1},\dots,\theta_{m-1},\sigma)\in 
\II_n^{m-t-1}\times S_n}\sgn(\sigma)\prod_{l=1}^n 
m_{k_{tl}\theta_{t+1}(l)}^{\tau^{m-1-t}}\dots 
m_{\theta_{m-1}(l)u_{\sigma(l)}}=0.$$}

\medskip

{\it Proof.} Assume that $k_{ti}=k_{tj}$ for $i\neq j$; consider the 
disjoint union
$$\II_n^{m-t-1}\times S_n=\II_n^{m-t-1}\times 
A_n\coprod\II_n^{m-t-1}\times S_n\backslash A_n.$$
The map 
$(\theta_{t+1},\dots,\theta_{m-1},\sigma)\mapsto(\theta_{t+1}\circ(i,j),\dots,\theta_{m-1}\circ(i,j),\sigma\circ(i,j))$

is a bijection from $\II_n^{m-t-1}\times A_n$ to $\II_n^{m-t-1}\times 
S_n\backslash A_n$. Since $k_{ti}=k_{tj}$ and $\sgn(\sigma)=-\sgn(\sigma\circ(i,j))$, the terms in $S''$ coming from 
$(\theta_{t+1},\dots,\theta_{m-1},\sigma)$ and 
$(\theta_{t+1}\circ(i,j),\dots,\theta_{m-1}\circ(i,j),\sigma\circ(i,j))$ 
cancel each other and we are done.

\medskip

Summing up, we get a new expression for $S_U$
$$S_U=\sum_{(\theta_1,\dots,\theta_{m-1})\in \II_n^{m-1}}\sum_{\sigma\in 
S_n}\sgn(\sigma)\prod_{i=1}^n 
m_{u_i\theta_1(i)}^{\tau^{m-1}}m_{\theta_1(i)\theta_2(i)}^{\tau^{m-2}}\dots 
m_{\theta_{m-1}(i)u_{\sigma(i)}}.$$

\medskip

We are ready to prove the following:

\medskip

{\bf Proposition 3.1} {\it Assume that $p\geq 3d$; then for any 
$1\leq n\leq d-1$, we have:
$$\M_n\equiv \sum_{(\sigma,\theta_1,\dots,\theta_{m-1})\in 
S_n^{m}}\sgn(\sigma)\prod_{i=1}^n 
m_{i\theta_1(i)}^{\tau^{m-1}}m_{\theta_1(i)\theta_2(i)}^{\tau^{m-2}}\dots 
m_{\theta_{m-1}(i)\sigma(i)}~[\pi^{mY_n+1}].$$}

\medskip

{\it Proof.} Let $V$ be the valuation of 
$m_{u_i\theta_1(i)}^{\tau^{m-1}}m_{\theta_1(i)\theta_2(i)}^{\tau^{m-2}}\dots 
m_{\theta_{m-1}(i)u_{\sigma(i)}}$; from Corollary 2.1 (note that since $d\geq 2$ and $p\geq 3d$ we have $p\geq d+3$), we get:
$$ V  \geq  \sum_{i=1}^n 
\lceil\frac{pu_i-\theta_1(i)}{d}\rceil+\dots+\lceil\frac{p\theta_{m-1}(i)-u_{\sigma(i)}}{d}\rceil$$
Assume that $1\leq u_1<\dots<u_n=n+t_0$, and $1\leq \theta_i(1)<\dots<\theta_i(n)=n+t_i$, $1\leq i\leq m-1$; then we have from lemma 2.3
$$\begin{array}{rcl}
V& \geq & Y_n+\left(\left[\frac{p}{d}\right]-1\right)t_0-t_1+\dots+y_n+\left(\left[\frac{p}{d}\right]-1\right)t_{m-1}-t_0\\
& \geq & mY_n+\left(\left[\frac{p}{d}\right]-2\right)(t_0+\dots+t_{m-1}).\\
\end{array}$$
Assume that one of the $t_i$ is nonzero; from the hypothesis on $p$, we have $V\geq mY_n+1$, and this term doesn't appear in the congruence. Thus the only terms remaining are 
those with $\{u_1,\dots,u_n\}$, $\theta_i(\{1,\dots,n\})$ all equal to 
$\{1,\dots,n\}$, and this is the desired result.

\medskip

We are now ready to show the main result of this section; we use the 
notations of section 2:

\medskip

{\bf Proposition 3.2} {\it Assume that $p\geq 3d$; then for any $1\leq 
n\leq d-1$, we have the congruence:
$$\M_n\equiv 
\frac{N_{\K_m/\ma{Q}_p}(\P_n(a_1,\dots,a_d))}{\left(\prod_{i\notin 
B_n}\lceil\frac{pi}{d}\rceil!\prod_{i\in 
B_n}\left(\lceil\frac{pi}{d}\rceil-1\right)!\right)^m}\pi^{mY_n}~[\pi^{mY_n+1}].$$}

\medskip

{\it Proof.} We rewrite the sum in proposition 3.1: set 
$\sigma_0=\theta_1$, 
$\sigma_1=\theta_2\circ\theta_1^{-1},\dots,\sigma_{m-1}=\sigma\circ 
\theta_{m-1}^{-1}$; we get

$$\M_n  \equiv  \sum_{(\sigma_0,\dots,\sigma_{m-1})\in 
S_n^m}\sgn(\sigma_0\circ\dots\circ\sigma_{m-1})\prod_{i=1}^n 
m_{i\sigma_0(i)}^{\tau^{m-1}}m_{i\sigma_1(i)}^{\tau^{m-2}}\dots 
m_{i\sigma_{m-1}(i)}~[\pi^{mY_n+1}]$$
$$\equiv \prod_{i=0}^{m-1} \sum_{\sigma_i\in 
S_n}\sgn(\sigma_i)\prod_{j=1}^n 
m_{j\sigma_i(j)}^{\tau^{m-1-i}}\equiv \prod_{i=0}^{m-1} \left(\sum_{\sigma_i\in 
S_n}\sgn(\sigma_i)\prod_{j=1}^n 
m_{j\sigma_i(j)}\right)^{\tau^{m-1-i}}~[\pi^{mY_n+1}].$$

Finally we know from Proposition 2.2 that
$$\begin{array}{rcl}
\left(\sum_{\sigma\in \Sigma_n}\sgn(\sigma)\prod_{j=1}^n 
m_{j\sigma(j)}\right)^{\tau^i} & \equiv & 
\left(\sum_{\sigma\in \Sigma_n}\sgn(\sigma)\prod_{j=1}^n 
m_{j\sigma(j)}\right)^{\tau^i} ~[\pi^{Y_n+1}]\\
& \equiv &
\frac{\P_n(a_1,\dots,a_d)^{\tau^i}}{\prod_{i\notin 
B_n}\lceil\frac{pi}{d}\rceil!\prod_{i\in 
B_n}\left(\lceil\frac{pi}{d}\rceil-1\right)!}\pi^{Y_n}~[\pi^{Y_n+1}],
\end{array}$$
and the theorem is an immediate consequence of the congruences above.

\medskip

\section{Generic Newton polygons}

In this section we use the results above to determine the generic Newton 
polygon $GNP(d,q)$ associated to polynomials of degree $d$ over 
$\F_q$. We determine the Zariski dense open subset $U$ in 
$\A^{d-1}$, the space of monic polynomials of degree $d$ without 
constant coefficient, such that for any $f\in U$ we have 
$NP_q(f,\F_q)=GNP(d,q)$, giving an explicit polynomial, the 
Hasse polynomial $G_{d,p}$ in $\F_p[X_1,\dots,X_d]$ such that 
$U=D(G_{d,p})$.
 
\medskip

\subsection{Hasse polynomials} In this section, we study the polynomials which 
appear when expressing the principal parts of the minors $M_n$ in terms 
of the coefficients of the original polynomial.

\medskip

{\bf Definition 4.1.} {\it Recall that for $g(X)=\sum_{i=1}^d a_iX^i$, 
we have set
$$\P_n(a_1,\dots,a_d):=\sum_{\sigma\in \Sigma_n} \sgn(\sigma)
\prod_{i=1}^n\left\{g^{\lceil\frac{pi-\sigma(i)}{d}\rceil}\right\}_{pi-\sigma(i)}.$$
We denote by $P_n\in \F_p[X_1,\dots,X_d]$ the reduction modulo $p$ of 
$\P_n$, and let $P_{d,p}:=\prod_{i=1}^{[\frac{d}{2}]} P_i$.}

\medskip

Our next task is to ensure that the polynomial $P_{d,p}$ is non zero; in order to prove this, we consider the monomials in $P_{d,p}$ of minimal degree and exhibit one that appear (with non zero coefficient) exactly once when $\sigma$ describes $\Sigma_n$.

\medskip

{\bf Lemma 4.1} {\it For any $1\leq n\leq d-1$, we have $P_n\neq 0$ in 
$\F_p[X_1,\dots,X_d]$. Moreover this polynomial is homogeneous of degree 
$Y_n$.}

\medskip

{\it Proof.} The polynomial $(a_1,\dots,a_d)\mapsto 
\left\{g^{\lceil\frac{pi-\sigma(i)}{d}\rceil}\right\}_{pi-\sigma(i)}$ 
contains a unique monomial of maximal degree in $X_d$, which is 
$X_d^{\left[\frac{pi-\sigma(i)}{d}\right]}X_{\overline{pi-\sigma(i)}}$, 
where $\overline{n}$ stands for the least nonnegative integer congruent 
to $n$ modulo $d$, and we set $X_0=1$. Moreover its coefficient is $1$ 
if $\overline{pi-\sigma(i)}=0$, and $\lceil\frac{pi-\sigma(i)}{d}\rceil$ 
else: in any case it is non zero modulo $p$. Thus 
$\prod_{i=1}^n\left\{g^{\lceil\frac{pi-\sigma(i)}{d}\rceil}\right\}_{pi-\sigma(i)}$ 
contains a unique monomial of maximal degree in $X_d$ with nonzero 
coefficient, which is 
$X_d^{\sum_{i=1}^n\left[\frac{pi-\sigma(i)}{d}\right]}\prod_{i=1}^n 
X_{\overline{pi-\sigma(i)}}$.

\medskip

On the other hand, we have 
$\left[\frac{pi-j}{d}\right]=\left[\frac{pi}{d}\right]$ if $1\leq j\leq 
j_i$, and $\left[\frac{pi-j}{d}\right]=\left[\frac{pi}{d}\right]-1$ if 
$j\geq j_i+1$. Thus the degree in $X_d$ of a monomial of $P_n$ is 
maximal for those $\sigma$ such that $\sigma(i)\leq j_i$. From Lemma 2.2,
 we see that the monomials in $P_n$ of maximal 
degree in $X_d$ come from the $\sigma$ such that for any $i\in B_n$, 
$\sigma(i)=j_i$ (note that such $\sigma$ exist since $i\mapsto j_i$ is 
injective on $B_n$). If $\Sigma_n^+\subset \Sigma_n$ is the set of these 
permutations, we get that the monomials in $P_n$ of maximal degree in 
$X_d$ are the
$$X_d^{\sum_{i=1}^n\left[\frac{pi}{d}\right]}\prod_{i\notin B_n} 
X_{\overline{pi-\sigma(i)}}=X_d^{\sum_{i=1}^n\left[\frac{pi}{d}\right]}\prod_{i\notin B_n} 
X_{j_i-\sigma(i)},$$
with $\sigma\in \Sigma_n^+$, and that there is exactly $\#\Sigma_n^+$ 
such monomials in $P_n$ (remark that for $i\notin B_n$, 
$\sigma\in \Sigma_n^+$, we have $\overline{pi}=j_i>n$, and 
$\overline{pi-\sigma(i)}=j_i-\sigma(i)$).

\medskip

We now construct $\sigma_0\in \Sigma_n^+$ such that the associated 
monomial cannot be obtained from any other $\sigma\in \Sigma_n^+$. For 
$i\in B_n$, we must have $\sigma_0(i)=j_i$ from the definition of 
$\Sigma_n^+$. Let $i_0 \in \{1,\dots,n\}\backslash B_n$ be such that 
$j_{i_0}$ is maximal, and set 
$\sigma_0(i_0)=\min\left\{\{1,\dots,n\}\backslash \{j_i,~i\in 
B_n\}\right\}$. Then we continue the same process, with $i_1\neq i_0$, 
$i_1\notin B_n$ such that $j_{i_1}$ is maximal, and $\sigma_0(i_1)$ 
minimal among the remaining possible images. The permutation $\sigma_0$ is 
clearly well defined, and unique. Let $\sigma\in \Sigma_n^+$ be 
such that $\prod_{i\notin B_n} X_{j_i-\sigma(i)}=\prod_{i\notin B_n} 
X_{j_i-\sigma_0(i)}$. Consequently here exists $i\notin B_n$ such that 
$j_i-\sigma(i)=j_{i_0}-\sigma_0(i_0)$; from the construction we must 
have $j_i=j_{i_0}$, thus $i=i_0$, and $\sigma(i_0)=\sigma_0(i_0)$. 
Following this process, we get $\sigma=\sigma_0$. Finally the monomial 
$X_d^{\sum_{i=1}^n\left[\frac{pi}{d}\right]}\prod_{i\notin B_n} 
X_{j_i-\sigma_0(i)}$ appears just once in $P_n$, with coefficient $\prod_{i\notin B_n}\lceil\frac{pi-\sigma_0(i)}{d}\rceil$ and this 
gives the first assertion.

\medskip

To prove the second assertion, remark that from the proof of Lemma 2.1, 
$(a_1,\dots,a_d)\mapsto 
\left\{g^{\lceil\frac{pi-\sigma(i)}{d}\rceil}\right\}_{pi-\sigma(i)}$ is 
homogeneous of degree $\lceil\frac{pi-\sigma(i)}{d}\rceil$; thus from 
the definition of $\Sigma_n$, we get the result.

\medskip

{\bf Lemma 4.2} {\it i) We have $P_{d,p}(X_1,\dots,X_{d-1},1)\neq 0$ in 
$\F_p[X_1,\dots,X_{d-1}]$. Moreover this polynomial has degree less or 
equal than $\frac{d-1}{2}\left[\frac{d}{2}\right]\left(\left[\frac{d}{2}\right]+1\right)$;

ii)  we have $P_{d,p}(X_1,\dots,X_{d-2},0,1)\neq 0$ in 
$\F_p[X_1,\dots,X_{d-1}]$. Moreover this polynomial has degree less or 
equal than $\frac{d-1}{4}\left[\frac{d}{2}\right]\left(\left[\frac{d}{2}\right]+1\right)$.}

\medskip

{\it Proof.} {\it i)} The non vanishing is obvious from Lemma 4.1, since 
dehomogeneizing a non zero homogeneous polynomial with respect to any 
of its variables yields a non zero polynomial. 

\medskip

We now show the assertion 
on the degree; consider the polynomial $(a_1,\dots,a_d)\mapsto 
\left\{g^{\lceil\frac{k}{d}\rceil}\right\}_{k}$. From the proof of Lemma 
2.1, its monomials are among the $X_1^{m_1}\dots X_d^{m_d}$ with 
$m_1+\dots+dm_d=k$, and $m_1+\dots+m_d=\lceil\frac{k}{d}\rceil$. 
Multiplying the second equality by $d$ and substracting the first we get 
$(d-1)m_1+\dots+m_{d-1}=d\lceil\frac{k}{d}\rceil-k\leq d-1$; 
consequently $m_1+\dots+m_{d-1}\leq d-1$, and the degree in 
$X_1,\dots,X_{d-1}$ of the above polynomial is at most $d-1$. From the 
definition of $\P_n$, its degree in the first $d-1$ variables is at most 
$n(d-1)$, and finally the degree of $P_{d,p}(X_1,\dots,X_{d-1},1)$ is at 
most $\frac{d-1}{2}\left[\frac{d}{2}\right]\left(\left[\frac{d}{2}\right]+1\right)$.

\medskip

{\it ii)} The non vanishing follows from the proof of Lemma 4.1. Remark that from the construction of $\sigma_0$, we must have, for $i\notin B_n$, $j_i-\sigma_0(i)\leq d-2$; thus the monomial constructed in the proof doesn't contain $X_{d-1}$, and the result follows. In order to give a bound for the degree, we use the same technique that in the proof of {\it i)}, remarking that now we take $m_{d-1}=0$, and consequently $m_1+\dots+m_{d-1}\leq \frac{d-1}{2}$. This ends the proof.

\medskip

{\bf Definition 4.2.} We define the {\it Hasse polynomial for polynomials of degree $d$} $G_{d,p}$ 
in $\F_q[X_1,\dots,X_{d-1}]$ as
$$G_{d,p}(X_1,\dots,X_{d-1}):=P_{d,p}(X_1,\dots,X_{d-1},1),$$
and the {\it Hasse polynomial for normalized polynomials of degree $d$}, $H_{d,p}$ 
in $\F_q[X_1,\dots,X_{d-2}]$ as
$$H_{d,p}(X_1,\dots,X_{d-2}):=P_{d,p}(X_1,\dots,X_{d-2},0,1),$$

\medskip

\subsection{The generic Newton polygon.} We use the results of the 
paragraph above to show that for any monic polynomial of degree $d$ over 
$\F_q$, its Newton polygon is above a generic Newton polygon, and that 
most polynomials have their Newton polygon attaining the generic Newton 
polygon.

\medskip

We identify the set of normalized monic polynomials of degree $d$ such that 
$f(0)=0$  with affine $d-2$ space $\A^{d-2}$ by associating the point 
$(a_1,\dots,a_{d-2})$ to the polynomial 
$f(X)=X^d+a_{d-2}X^{d-2}+\dots+a_1X$.
 
\medskip

{\bf Definition 4.3.} Set $Y_0:=0$. We define the {\it generic Newton 
polygon} of exponential sums associated to polynomials of degree $d$ in 
$\F_q$, $GNP(d,\F_q)$, as the lowest convex hull of the points
$$\left\{ (n,\frac{Y_n}{p-1})\right\}_{0\leq n \leq d-1}.$$

\medskip

We are ready to prove the main result of this paper.

\medskip
 
{\bf Theorem 4.1.} {\it Let $p\geq 3d$ be a prime, and $f\in \F_q[X]$ 
a normalized polynomial of degree $d$. Then we have $NP_q(f,\F_q)=GNP(d,q)$ if and only if the coefficients of $f$ belong to the Zariski dense open subset $U:=D(H_{p,d})$. 
Moreover for any polynomial of degree $d$ over $\F_q$, the associated 
Newton polygon is above the generic Newton polygon.}

\medskip

{\it Proof.} Recall from Proposition $1.1$ that for any polynomial of degree $d$ we have
$$L(f,T)=\det(\I_{d-1}-T\Gamma^{\tau^{m-1}}\dots \Gamma)=\sum_{n=0}^{d-1} \M_nT^n.$$
Thus the Newton polygon $NP_q(f,\F_q)$ is the lower convex hull of the set of points 
$$\{ (n,v_q(\M_n)),~0\leq n\leq d-1\}.$$
On the other hand we have, from Proposition 3.2 
$$\M_n\equiv 
\frac{N_{\K_m/\ma{Q}_p}(\P_n(a_1,\dots,a_d))}{\left(\prod_{i\notin 
B_n}\lceil\frac{pi}{d}\rceil!\prod_{i\in 
B_n}\left(\lceil\frac{pi}{d}\rceil-1\right)!\right)^m}\pi^{mY_n}~[\pi^{mY_n+1}].$$
and we get $v_q(\M_n)=\frac{Y_n}{p-1}$ if and only if $P_n(\alpha_1,\dots,\alpha_{d-2},0,1)\neq 0$ in $\F_q$. 
Moreover, the Newton polygon is symmetric : if it has a slope of length $l$ and slope $s$ it has a segment of the same length and slope $1-s$. Thus, in order to show that $NP_q(f,\F_q)$ coincides with $GNP(d,q)$, it is sufficient to show that the first $[\frac{d}{2}]$ vertices of $NP_q(f,\F_q)$ coincide with the ones of $GNP(d,q)$. From Definition 4.1, this is true exactly when $P_{d,p}(\alpha_1,\dots,\alpha_{d-2},0,1)\neq 0$; this is the desired result. The last assertion is an easy consequence of the discussion above.

\bigskip

{\bf Remark 4.1.} Let us show that we have $NP(f)=HP(d)$ for any $f$ of degree $d$ when $p\equiv 1~[d]$; in this case we get $\Sigma_n=\{Id\}$ for any $n$, $Y_n=(p-1)\frac{n(n+1)}{2}$ and $GNP(d,q)=HP(d)$; moreover $P_n(X_1,\dots,X_d)=cX_d^{Y_n}$ for some $c\in \F_p^\times$, and $H_{d,p}$ is a nonzero polynomial of degree $0$. In this case we get that $U_{d,p}$ is the whole $\A^{d-2}$, as stated above.

\end{document}